\documentclass{article}
\usepackage{amsmath}
\usepackage{amssymb}
\usepackage{amsfonts}
\usepackage{theorem}
\theorembodyfont{\upshape}

\newtheorem{thr}{Theorem}[section]

\newtheorem{lm}{Lemma}[section]

\newtheorem{cor}{Corollary}[section]
\usepackage{verbatim}

\newcommand{\ri}{\rightarrow}
\newcommand{\eps}{\varepsilon}
\newcommand{\E}{\mathsf{E}}
\renewcommand{\P}{\mathsf{P}}
\newcommand{\R}{{\textbf{R}}}
\newcommand{\Z}{{\textbf{Z}}}
\newcommand{\F}{{\cal F}}

\newcommand{\by}{{\bf y}}
\newcommand{\bx}{{\bf x}}

\date{}

\begin{document}

\title{{\bf Limit theorems for  random point measures 
generated by cooperative sequential adsorption}}

\author{
V. Shcherbakov,\\
{\small Laboratory of Large Random Systems},\\
{\small  Faculty 
of Mathematics and Mechanics,}\\ 
{\small Moscow
State University, 119992, Moscow, Russia,}\\ 
{\small Shcherbakov@mech.math.msu.su}}

\maketitle

\vspace*{-5ex}
\begin{abstract}
We consider a finite sequence of random points 
in a finite domain of a finite-dimensional Euclidean space.
The points are sequentially allocated in the domain 
according to a model of cooperative sequential adsorption.
The main peculiarity of the model is that  
the probability distribution of a point depends on previously 
allocated points. We assume that the dependence vanishes
as the concentration of points tends to infinity. Under this assumption 
the law of large numbers, the central limit theorem 
and Poisson approximation are proved
for the generated sequence of random point measures.

   \medskip
  {\bf 2000 MSC}: 82C21 (Primary) 60F05, 60F15, 60G55 (Secondary).

%  \medskip
%  {\it Keywords}: cooperative sequential adsorption;
%infinite range cooperative effects;
%the law of large numbers; Poisson approximation; 
%the central limit theorem; Gaussian random field.

\end{abstract}

\section{Introduction and the results}
\label{dynamo}
In this paper we study the asymptotic behavior 
of random point measures 
\begin{equation}
\label{mes1}
\mu_m=\sum\limits_{i=1}^{m}\delta_{X_i},
\end{equation}
generated  by random points $X_1,\ldots,X_m$ sequentially allocated 
in a compact set $D\subset \R^d$.
To describe the joint distribution of  $X_1,\ldots,X_m$
we need some notation.
For any point $x\in D$ and a finite non-empty 
set $\by=\{y_1,\ldots,y_n\},\, n\geq 1,$ of points in $D$ we 
denote by $n(x,\by)$ the number of points $y_i\in \by$,
such that the distance between $x$ and $y_i$ is not greater than $R(x)$, 
where $R:D\ri \R_{+}$ is some  measurable function.
By definition $n(x,\emptyset)=0.$
The number $R(x)$ is called an interaction radius at point $x$.
Let $\{\beta_{n}(x),\, n\geq 0\}$ be  a sequence 
of measurable positive bounded functions on $D$.
 Denote for short $X(k)=(X_1,\ldots,X_k),\, k\geq 1,$ and $X(0)=\emptyset$.
Given the set of points $X(k)$ the conditional 
distribution of point $X_{k+1}$ is specified by the 
following probability density
\begin{equation}
\label{f1}
\psi_{k+1}(x)=\frac{\beta_{n(x,X(k))}(x)}{\alpha(X(k))},
\end{equation} 
where
$$\alpha(X(k))=\int\limits_{D}\beta_{n(u,X(k))}(u)du,$$ 
is the normalizing constant.
The joint probability density  of $X_1,\ldots,X_m$ at points
$x_1,\ldots,x_m$ is 
\begin{equation}
\label{cond}
p_m(x_1,\ldots,x_m)=\prod_{k=1}^{m}\frac{\beta_{n(x_k, \bx_{<k})}(x_k)}
{\int_{D}\beta_{n(x, \bx_{<k})}(x)dx}
=\prod_{k=1}^{m}\psi(x_k|\bx_{<k}),
\end{equation}
where we denoted for short 
 $\bx_{<k}=(x_1,\ldots,x_{k-1}),\, k\geq 2,$
and $\bx_{<1}=\emptyset$ for $k=1$.

Let us give examples of the situation where  this  set of 
sequentially allocated random points naturally appears.
First we do it  in terms of  dynamic processes with continuous
time  describing adsorption 
reactions with cooperative effects. Namely,
consider a spatial
birth process $\bx(t),\,t\geq 0,$ in $D$ with
birth rates defined in terms of functions 
$\beta_{n}(x),\,n\geq 0,$ as follows. If the process state 
at time $t\geq 0$ is $\bx$, then  the birth rates are 
$\beta_{n(x,\bx)}(x),\, x\in D$, so
the total birth rate is  $\alpha(\bx)$ 
and the time until the next jump is an 
exponential random variable with mean  $\alpha^{-1}(\bx)$.
Assume that $\bx(0)=\emptyset$ and consider a random point process 
 $X(m)=(X_k,\, k=1,\ldots,m)$ formed by the first $m$ points
of the spatial birth process $\bx(t)$. 
It is easy to see that the first point $X_1$ has the probability
distribution specified by the function $\beta_0(x)$
normalized to be a probability density.
Given  $X_1,\ldots,X_k,\, k\geq 1,$
the conditional distribution of $X_{k+1}$ is specified 
by the probability density  (\ref{f1}).
The spatial birth process just described 
is a continuous version of  a lattice model 
of monomer filling with nearest-neighbor cooperative effects (\cite{Evans}).
It is a particular case of the models of cooperative
sequential adsorption widely used in 
physics and chemistry for modeling various adsorption processes
(see \cite{Evans} and \cite{Privman} for more details and 
surveys of the relevant literature). 
%(see refs. $2$ and $5$ for more details and 
%surveys of the relevant literature).
The set of random points $X(m)$ can also  be viewed as an output
of the following sequential packing process with discrete time.
Consider a sequence of random points $Y_i,\, i\geq 1,$
sequentially  arriving in $D$. Each point $Y_i$ is uniformly distributed
in $D$ and is accepted with probability depending
on a number of previously accepted points in the local configuration
near $Y_i$. More precisely, let $Y(N)=(Y_1,\ldots, Y_N)$ be a set of the 
first $N$ arrived points and let $X(k)=(X_1,\ldots, X_k),\, k=k(N),$ be a set
of accepted ones among $Y_i,\, i=1,\ldots,N.$ Next uniformly distributed
arrival $Y_{N+1}$ is accepted with probability   
$\beta_{n(Y_{N+1}, X(k))}(Y_{N+1})/C,$ where $C$ is an arbitrary constant
such that $\sup_{n}\sup_{x\in D}\beta_n(x)\leq C$. 
Regardless of a particular choice of $C$ 
the probability density  of the next {\it accepted} point $X_{k+1}$
is given by the formula (\ref{f1}).
The value of $C$ influences only a number of discarded arrivals $Y'$s
until next acceptance.
Thus, given  the set of previously accepted points $X(k)$, we 
use  a well known acceptance-rejection 
sampling  for simulating  a random variable  
which distribution is specified by the unnormalized
probability density $\beta_{n(x, X(k))}(x),\, x\in D$.
The sequence  of points $X(m)$ is a set of first $m$ 
sequentially accepted points.   
  
The measures (\ref{mes1})
belong to the class  
of random point measures generated  by the spatial processes
arising in random sequential packing and deposition problems
(see \cite{BarYuk}, \cite{PenYuk} and references therein).
The typical  example is when  one
sequentially allocates  $m$ points in a unit cube.
Each point is uniformly distributed in the cube and is
accepted with probability depending on configuration
of previously accepted points in  the ball of radius $1/m$ around the
point.  
Therefore,  the interaction radius in those models is inversely 
proportional to the number of points and this leads to  
the well-known effect of  finite range dependence between points.
It is not the case in our model where 
the interaction radius is a fixed positive function (or constant) 
regardless of the number of points. This corresponds
to the so-called  infinite range of interaction or infinite range 
 cooperative effects, see, for instance, \cite{Evans}. 

Our other main assumption is that 
$\beta_n(x) \ri \beta(x)>0$ as $n\ri \infty$ uniformly in $x\in D$,
where the function $\beta$ is bounded below and above.
Under our assumptions the sequence of random variables
$X_k,\,k\geq 1,$ converges in total variation to a random variable
with the probability density specified by  the function $\beta(x),\,x\in D,$
appropriately normalized.  
Therefore the model can be considered as a perturbation 
of the binomial case which is
$\beta_n(x)=\beta(x),\, x\in D$ for any $n\geq 0$.
The perturbation vanishes while the domain is saturated by points.
The distribution of  a new arrival   becomes "more uniform" and
"more independent" on the existing configuration of points provided 
the domain is sufficiently saturated and the saturation
is "sufficiently uniform". We make it rigorous in Lemma \ref{L1}. 
From the physical point of view the assumption on the sequence 
of intensities can be interpreted as follows.
One might think of an adsorption process 
such that reaction rates depend on local environment and 
stabilize  when the concentration of adsorbed molecules is sufficiently high.

In the binomial case we immediately get 
Theorems \ref{M1}, \ref{Pois} and \ref{M2},
since the points are independent.
In general case the points  are  dependent
and we arrive at the proof of the law of large numbers,
the central limit theorem and Poisson approximation for the sequence 
of {\it dependent} random variables.
Some care should be taken 
to assess  the weakening of dependence in the tail of the sequence
$X(m)$. 
Note that we obtain  the central limit theorem (Theorem \ref{M2})
assuming that the 
sequence of functions $\{\beta_n(x),\, n\geq 0\}$ converges
to its limit with some rate. 

{\bf Remark.} We will denote by the  letter $C$ or by the letter $C$ 
with subscripts the
various constants the  particular values of which are immaterial 
for the proofs. In some cases we will stress dependence 
of these constants on some parameters that do not depend 
on the number of points $m$.   
By ${\cal B}(D)$ the set of real-valued 
measurable bounded functions on $D$ is denoted  and $\|f\|_{\infty}=
\sup_{x\in D}|f(x)|$ for $f\in {\cal B}(D)$.
It is assumed that the random variables $X_k,\, k\geq 1,$ are realized on some 
probability space with probability measure  $\P$ and 
$\E$ is expectation with respect to $\P$.

\begin{thr}
\label{M1}
Assume that $\inf_{x\in D}R(x)>0$, the  sequence of positive
functions $\beta_n\in {\cal B}(D),\,n\geq 0,$
is uniformly bounded and converges uniformly as $n\ri \infty$  
to a function   $\beta \in {\cal B}(D)$, such that 
$\inf_{x\in D}\beta(x)>0.$
Then the law of large numbers holds for the sequence of random
measures $\mu_m$. That is  
for any function $f\in {\cal B}(D)$
$$\frac{1}{m}\int\limits_{D}f(x)d\mu_{m}(x)=
\frac{1}{m}\sum\limits_{i=1}^{m}f(X_i) \ri  
J(f)=\frac{1}{\alpha}\int\limits_{D}f(x)\beta(x)dx,$$
in probability as $m \ri \infty$,
where  $\alpha=\int_{D}\beta(x)dx.$
\end{thr}
\begin{thr}
\label{Pois}
In addition to the assumptions of Theorem \ref{M1} assume that
the function $\beta$ is continuous. Fix an arbitrary
$x\in D$ and $r>0$. Let $S_{m}(x,r)$ be a number
of those points $X_k,\, k=1,\ldots,m,$ that fall
in a ball $B(x,rm^{-1/d})$. Then 
a sequence of random variables $S_{m}(x,r),\, m\geq 1,$
converges in a weak sense to a Poisson random
variable with parameter $r^db_d\beta(x)/\alpha$,
where $b_d$ is a volume of a $d-$dimensional ball with unit radius. 
\end{thr}

\begin{thr}
\label{M2}
In addition to the assumptions of Theorem \ref{M1} assume
that 
\begin{equation}
\label{expon}
|\beta_n(x) - \beta(x)|\leq \tau(x)\varphi(n),
\end{equation} 
for any $n\geq 0$, where a function
 $\varphi(s)>0,\, s\geq 0,$ is such that $\varphi(s)\ri 0$ as 
$s\ri \infty$ and  for any $\delta>0$
\begin{equation}
\label{var}
\frac{1}{\sqrt{n}}\sum\limits_{k=0}^{n}\varphi(k\delta) \ri 0,
\end{equation}
as $n\ri \infty$,
the function $\tau \in {\cal B}(D)$ is such that
$\inf_{x\in D} \tau(x)>0$.
Then the sequence of centred and rescaled random measures 
$(\mu_m-\E \mu_m)/\sqrt{m}$
converges as $m \ri \infty$ to a generalized Gaussian random field
on $D$ with zero mean and the covariance kernel  
\begin{align*}
G(f,g)&=J(fg)-J(f)J(g)\\
&=
\frac{1}{\alpha}\int\limits_{D}f(x)g(x)\beta(x)dx - 
\frac{1}{\alpha^2}\int\limits_{D}f(x)\beta(x)dx\int\limits_{D}g(x)\beta(x)dx,
\end{align*}
for any functions $f, g \in {\cal B}(D)$.
\end{thr}
To prove these theorems we will use Lemmas \ref{L1}-\ref{L4}.
\begin{lm}
\label{L1}
Assume that $\inf_{x\in D}R(x)>0$ and
$$0<\beta_{min}=\inf_{n}\inf\limits_{x\in D}
\beta_n(x) \leq \beta_{max}=\sup_{n}\sup\limits_{x\in D}\beta_n(x) < \infty,$$
then there exists a positive
constant $\delta_0$ such that 
for any  $\delta\in (0, \delta_0)$ 
\begin{equation}
\label{neq0}
\P\left\{\inf\limits_{x\in D} n(x,X(m))\leq m\delta\right\}
\leq Ce^{-\lambda m},
\end{equation} 
with some positive constants $C=C(\delta)$ 
and $\lambda=\lambda(\delta)$ 
for all sufficiently large $m$.

If the  assumptions of Theorem \ref{M1} hold, then for 
any $\eps>0$  
\begin{equation}
\label{neq1}
\P\left\{\sup_{x\in D} |\beta_{n(x,X(m))}(x) - \beta(x)|\geq \eps \right\}
\leq Ce^{-\lambda m},
\end{equation} 
and
\begin{equation}
\label{neq2}
\P\left\{|\alpha (X(m)) - \alpha|\geq \eps \right\}\leq Ce^{-\lambda m}
\end{equation} 
with the same  positive constants $C$ 
and $\lambda$ for all sufficiently large $m$.
\end{lm}
\begin{cor}
\label{c1}
If the  assumptions of Theorem \ref{M1} hold, then  
the sequence $X_m,\,m\geq 1,$ converges in total variation to a 
random variable $X$ distributed
according to the density $\beta(x)/\alpha$, as $m \ri \infty$. 
\end{cor} 
Let ${\cal F}_{k-1}$ be  a $\sigma-$algebra generated by the random variables 
$X_1,\ldots,X_{k-1}.$
For any function $f\in {\cal B}(D)$ denote 
$$J_k(f)=\E(f(X_k)|{\cal F}_{k-1}).$$

\begin{lm}
\label{L2}

1) If the assumptions of Theorem \ref{M1} hold, then for any 
 function $f \in {\cal B}(D)$ and for any  $p\geq 1$ 
 $$\E |J_k(f)- J(f)|^p \ri 0$$
as $k \ri \infty.$

2) If the  assumptions of Theorem \ref{M2} hold
and $\delta_0$ is the constant determined in Lemma \ref{L1}, 
then for any  $\delta\in(0, \delta_0)$
$$\E |J_k(f)- J(f)|^p \leq C \left(\varphi^p(k\delta)+
e^{-\lambda k}\right)
$$
as $k \ri \infty,$ with some constant $\lambda=\lambda(\delta)$.
\end{lm}
Let $Y$ be a random variable with probability density
$\beta(x)/\alpha$.
For any  function $f\in {\cal B}(D)$ and $n\geq 1$ denote
\begin{equation}
\label{fc}
{\cal U}_n(f)=\E(f(Y)-\E f(Y))^n=
\sum\limits_{i=0}^{n}(-1)^{n-i}{n\choose i} J(f^i)J^{n-i}(f),
\end{equation}
and   $\xi_k(f)=f(X_k)-\E f(X_k)$.
\begin{cor}
\label{c2}
Let $f\in {\cal B}(D)$ and fix some positive integer $n$. Then 

1) under assumptions of Theorem \ref{M1}
$$\E \left|\E\left(\xi_k^n(f)|{\cal F}_{k-1}\right)
- {\cal U}_n(f)\right| \ri 0 $$
as $k\ri \infty$, and 

2) under assumptions of Theorem \ref{M2}
$$\E \left|\E\left(\xi_k^n(f)|{\cal F}_{k-1}\right)
- {\cal U}_n(f)\right| \leq C \left(\varphi(k\delta)+e^{-\lambda k}\right)$$
as $k\ri \infty$. 
\end{cor} 

\begin{lm}
\label{L3}
Fix  a set of functions 
$g_{1},\ldots,g_{k} \in {\cal B}(D)$ and a set of positive integers 
 $r_1,\ldots,r_k$ and let $n=r_1+\cdots r_k$. Let  a set of
indices be such that $i_1<\cdots<i_k$ and denote by  $\eta$  
a random variable measurable with respect to the $\sigma-$algebra
${\cal F}_{i_1-1}$.

1) If the  assumptions of Theorem \ref{M1} hold, then
$$\left|\E \left(\eta \prod\limits_{v=1}^{k}\xi^{r_k}_{i_v}(g_{v})\right)-
 \E\eta \left(\prod\limits_{v=1}^{k}{\cal U}_{r_k}(g_{v})\right)
\right| \ri 0$$
as $i_1\ri \infty$. In particular, for any $f,g\in {\cal B}(D)$ and $k\neq j$
\begin{displaymath}
Cov(f(X_k),g(X_j))\ri 0,
\end{displaymath}
 as $\max(k,j)\ri \infty$. 

2) If the assumptions of Theorem \ref{M2} hold, then 
there exist constants $C=C(k,g_{1},\ldots,g_{k})$ 
such that for any $\delta \in (0,\delta_0)$
\begin{equation}
\label{bb}
\left|\E \eta \prod\limits_{v=1}^{k}\xi^{r_k}_{i_v}(g_{v})-
 \E \eta \prod\limits_{v=1}^{k}{\cal U}_{r_k}(g_{v})\right| 
\leq C^n\sum\limits_{v=1}^{k} 
\left(\varphi(i_v\delta)+e^{-\lambda i_v}\right),
\end{equation}
for all sufficiently large indices  
$i_1<\ldots<i_k,\, k\geq 1,$
where the constants 
$\lambda$ and $\delta_0$ are determined in Lemma \ref{L1}.
\end{lm}

\begin{lm}
\label{L4}
Under the assumptions of Theorem \ref{M2}
the sequence of random variables 
$$\frac{1}{\sqrt{m}}\sum\limits_{k=1}^{m}
(J_k(f)-\E f(X_k))$$
converges  to $0$ in probability as $m\ri \infty$.
\end{lm}

\section{Proofs}
\label{proflln}

{\bf Proof of Theorem \ref{M1}.}
Let us prove first that for any function $f\in {\cal B}(D)$ 
\begin{equation}
\label{eq1m1}
\frac{1}{m}\sum\limits_{k=1}^{m}\E f(X_k) \ri J(f),
\end{equation}
as $m\ri \infty$.
Indeed, By Lemma \ref{L2} we have
that $\E f(X_k) \ri J(f)$, as $k\ri \infty$.
Fix an arbitrary $\eps>0$ and let $k(\eps)$ be such that 
$|\E f(X_k)-J(f)|\leq \eps$ as $k>k(\eps)$.
It is easy to see that 
$$\left|\frac{1}{m}\sum\limits_{k=1}^{m}\E f(X_k) - J(f)\right|\leq
2\frac{k(\eps)}{m}\|f\|_{\infty}+\frac{m-k(\eps)}{m}\eps.$$
The first term in the right side of the preceding equation goes
to $0$ as $m\ri \infty$, the second is less than $\eps$.
Thus we get (\ref{eq1m1}) since $\eps$ is arbitrary.
It suffices now to prove that 
$$\frac{1}{m}\sum\limits_{k=1}^{m}(f(X_k)-\E f(X_k)) \ri 0,$$
in probability as $m\ri \infty$.
By Chebyshev inequality we have that for any $\eps>0$   
$$\P\left\{\left|\sum\limits_{k=1}^{m}f(X_k)-\E f(X_k)\right| \geq \eps m
\right\} \leq \frac{1}{\eps^2m^2} 
\sum\limits_{k,j=1}^{m}Cov(f(X_k),f(X_j)).$$
If $k\neq j$, then by part $1)$ of Lemma 
\ref{L3} $Cov(f(X_k), f(X_j))\ri 0$ as $\max(k,j) \ri 0$,
therefore the right hand side of the preceding 
display vanishes as $m\ri \infty$.
Theorem \ref{M1} is proved.

{\bf Proof of Theorem \ref{Pois}.}
Let $x\in D$ and $r>0$ be  fixed. Denote for short
$S_m=S_m(x,r)$. 
We prove that for any $t\in \R$
\begin{equation}
\label{ch0}
\lim\limits_{m\ri \infty}\E e^{itS_m}=\exp\{(e^{it}-1)\beta(x)r^db_d/\alpha\}.
\end{equation}
By definition
$$S_m=\sum\limits_{k=1}^{m}\xi_{m,k},$$
where $\xi_{m,k}=1_{\{X_k\in B(x,rm^{-1/d})\}}$.
For any $k\geq 1$ we can write
\begin{equation}
\label{ch}
\E \left(e^{it\xi_{m,k}}|{\cal F}_{k-1}\right)=1+
\left(e^{it} - 1\right)p_m
+\left(e^{it} - 1\right)(p_{m,k}-p_m),
\end{equation}
where $p_{m,k}=\P\{X_k\in B(x,rm^{-1/d})|{\cal F}_{k-1}\}$
and $p_m$ is the probability that a random variable
with density $\beta(y)/\alpha, y\in D,$ falls in the ball $B(x,rm^{-1/d})$.
Repeatedly using the equation (\ref{ch})  we obtain that 
\begin{multline*}
\E e^{itS_m}=\left(1+\left(e^{it}-1\right)p_m\right)^{m}
\\+\left(e^{it}-1\right)\sum\limits_{k=1}^{m}
\left(1+\left(e^{it}-1\right)p_m\right)^{m-k}(\E p_{m,k}-p_m).
\end{multline*}
It is easy to see that 
$mp_m\ri \beta(x)r^db_d/\alpha$
as $m\ri \infty$. Therefore the first term in the left hand  side of 
the preceding display
tends to the characteristic function of the Poisson distribution
with parameter $\beta(x)r^db_d/\alpha$.
Let us show that the second term in  the left hand  side of 
the preceding display vanishes as $m\ri \infty$.
Noting that 
$p_m=J(f_m)$
and $p_{m,k}=J_k(f_m)$ with function
$f_m(y)=1_{\{y\in B(x,rm^{-1/d})\}}$
and using Remark after the proof of Lemma \ref{L2} (the bound (\ref{eks00}))
we can write
\begin{equation}
\label{pmk}
\E |p_{m,k}-p_m |\leq
\frac{C}{m} \E \sup\limits_{y\in D}|\beta_{n(y,X(k-1))}(y)-\beta(y)|.
\end{equation}
Fix an arbitrary $\eps>0$. An argument leading
to the bounds (\ref{eq3}) and (\ref{eq4}) in the proof
of Lemma \ref{L2} gives us here  that  there exists such $k(\eps)$
that for any $k\geq k(\eps)$
 we can replace the bound (\ref{pmk}) by the following one
\begin{equation}
\label{pmk1}
 \E |p_{m,k}-p_m |\leq \frac{C_1}{m} \left(\eps + e^{-\lambda k}\right),
\end{equation}
where  constant $\lambda$ is the same as in Lemma \ref{L1}. Hence we can bound
$$\left|\left(e^{it}-1\right)\sum\limits_{k=1}^{m}
\left(1+\left(e^{it}-1\right)p_m\right)^{m-k}(\E p_{m,k}-p_m)\right|
\leq \left(C_2 \eps + \frac{C_3}{m}\right).$$
Therefore we finished the proof since $\eps$ was taken arbitrary.

{\bf Remark.} Using Theorem 1 in \cite{Serf} 
(a general result on Poisson approximation 
for sums of possibly dependent nonnegative  integer-valued
random variables) one can also bound 
\begin{equation}
\label{poiss}
\sup\limits_{A\subset \Z_{+}}
\left|\P\{S_m\in A\} - \P\{Y_m\in A\}\right| \leq \sum\limits_{k=1}^{m} p_{m}^{2}
+  \sum\limits_{k=1}^{m} \E |p_{m,k} - p_m|, 
\end{equation}
where $Y_m$ is a Poisson random variable with parameter 
$mp_m$. Combining the bound (\ref{pmk1}) with the fact that
$mp_m$ has a finite limit
as $m\ri \infty$ one can show that the right hand side of the equation
(\ref{poiss}) vanishes as $m\ri \infty$.

{\bf Proof of Theorem \ref{M2}.}
It suffices to prove that for any function $f\in {\cal B}(D)$
the sequence of random variables 
\begin{equation}
\label{smf}
S_m(f)=\frac{1}{\sqrt{m}}\sum\limits_{k=1}^{m}(f(X_k)-\E f(X_k))
\end{equation}
converges weakly  as $m\ri \infty$ to a Gaussian random variable 
with mean zero and the variance $G(f,f)=J(f^2)-J^2(f)$.
Note that 
\begin{equation}
\label{eq1m2}
S_m(f)=Z_m(f)+\frac{1}{\sqrt{m}}\sum\limits_{k=1}^{m}
(J_k(f)-\E f(X_k)),\, m\geq 1,
\end{equation}
where 
\begin{align*}
\label{zmf}
Z_m(f)& = \frac{1}{\sqrt{m}}\sum\limits_{k=1}^{m}
(f(X_k)-\E(f(X_k)|{\cal F}_{k-1}))\\
& =\frac{1}{\sqrt{m}}\sum\limits_{k=1}^{m}
(f(X_k)-J_k(f)),\, m\geq 1.
\end{align*}
By Lemma \ref{L4} the second term in the right hand side
of the equation (\ref{eq1m2}) converges to $0$ as $m\ri \infty$.
Therefore to prove the theorem we need to prove 
that the sequence of random variables 
$Z_m(f),\, m\geq 1,$ 
converges weakly to a Gaussian random variable with mean zero and 
the variance $G(f,f)$  as $m\ri \infty$.
Note that $\{Z_m(f),\, {\cal F}_{m},\, m\geq 1\}$
is a zero-mean, square-integrable martingale array with differences
$\zeta_{mk}=(f(X_k)-J_k(f))/\sqrt{m},\, k=1,\ldots,m$.
It is easy to see that
\begin{equation}
\label{i}
\max\limits_{k}|\zeta_{mk}|\leq \frac{2\|f\|_{\infty}}{\sqrt{m}}
\ri 0,
\end{equation}
and
\begin{equation}
\label{ii}
\E \left(\max\limits_{k}\zeta^2_{mk}\right)\leq \frac{4\|f\|^2_{\infty}}{m}
\ri 0.
\end{equation}
By Corollary \ref{c1} and Lemma \ref{L2}
$\E (f(X_k)-J_{k}(f))^2$ converges to $G(f,f)$
as $k\ri \infty$. Consequently 
$\sum_{k=1}^{m}\E \zeta_{mk}^2$ converges to $G(f,f)$
as $m\ri \infty$.
Combining the results of Lemmas \ref{L2}
and \ref{L3} it is easy to obtain that 
$Cov((f(X_k)-J_{k}(f))^2,(f(X_j)-J_{j}(f))^2$ tends to
$0$ for $k\neq j$ as $\max(k,j)\ri \infty$. It yields that
$Var\left(\sum_{k=1}^{m} \zeta^2_{mk}\right)$ vanishes 
as $m\ri \infty$. Therefore
\begin{equation}
\label{iii}
\sum\limits_{k=1}^m\zeta^2_{mk} \ri G(f,f),
\end{equation}
in probability as $m\ri \infty$.

The equations (\ref{i}), (\ref{ii})  and (\ref{iii})
mean that the conditions  of Theorem 3.2
in \cite{Hall}
hold for the martingale array  $\{Z_m(f),\, {\cal F}_{m},\, m\geq 1\}$.
Therefore $Z_m(f)$ converges in distribution to a
Gaussian random variable with zero mean and covariance
$G(f,f)$ as $m\ri \infty$ and Theorem \ref{M2} is proved.

{\bf Proof of Lemma \ref{L1}.}
Without loss of generality 
we assume 
that the set $D$ is a $d-$dimensional unit cube.
If $l\in \Z_{+}$ is  the minimal integer such that
$$p(l)=l^{-d}\frac{\beta_{min}}{\beta_{max}}<1,\,\,\, \mbox{and}\,\,\,
1/l < \frac{1}{4}\inf\limits_{x\in D}R(x),$$
then we put $\delta_0=p(l)$.
Let $\{Q_i,\, i=1,\ldots, l^d\}$ be a set of 
non-overlapping cubes of size $1/l$ such that
$D=\bigcup_{i}Q_i.$  
Denote by $\xi_{mi}$ a number of points $X_1,\ldots,X_m$
falling in the cube $Q_i$.
Take a point $x\in D$ and let $x\in Q_i$ for some $i$.
It is easy to see that 
\begin{equation}
\label{eq1}
n(x,X(m))\geq \xi_{mi} \geq \min\limits_{j}\xi_{mj},
\end{equation}
since $Q_{i} \subset B(x,R(x))$.
The equation (\ref{eq1}) implies that 
$$\{n(x, X(m))\leq z\}\subset A_m=
\left\{\min\limits_{i}\xi_{mi}\leq z\right\}=
\bigcup\limits_{i}\{\xi_{mi}\leq z\},$$
   for any $z>0$.
It is obvious that 
$$\P\{A_{m}\}\leq l^d \max\limits_{i}
\P\{\xi_{mi} \leq z\}.$$
The formula (\ref{f1}) yields that
$$\P\{X_k\in Q_i| X(k-1)\}=\frac{\int_{Q_i}\beta_{n(u,X(k-1))}(u)du}
{\int_{D}\beta_{n(u,X(k-1))}(u)du}.$$
This conditional probability can be bounded  
below by $p(l)$ uniformly in sequences $X(k-1)$. Therefore 
the  unconditional probability $\P\{X_k\in Q_i\}$ 
is also bounded below  by the same constant for any $k\geq 1$.
Using the well-known coupling construction
%(e.g. Theorem~I.5.1 in ref. $8$)
we can construct on the same probability space 
the random variable $\xi_{mi}$ and 
 the binomial random variable $\tilde{\xi}_{mi}$  with $m$ trials and 
with  $p(l)$ the probability of success such that 
$\xi_{mi}$ stochastically dominates
$\tilde{\xi}_{mi}$.
So, we have that 
\begin{displaymath}
\P\{\xi_{mi}\leq m\delta\} \leq \P\{\tilde{\xi}_{mi}\leq m\delta\}
\end{displaymath}
for any $\delta>0$. If we take $\delta$ such that $0<\delta<\delta_0=p(l)$, then 
the well known  large deviations bounds 
%(see, for instance, ref. $3$)
for  the sums of i.i.d. random variables give us that 
$$\P\{\tilde{\xi}_{mi}\leq m\delta\} \leq Ce^{-\lambda m},$$
with some positive constants $C$ and $\lambda$.
Therefore
\begin{equation*}
\P\left\{\inf\limits_{x\in D} n(x,X(m))\leq m\delta\right\}\leq 
l^d\max\limits_{i}\P\{\xi_{mi}\leq m\delta\} \leq 
Cl^de^{-\lambda m}
\end{equation*}
and the proof of the bound (\ref{neq0}) is over.
The bounds (\ref{neq1}) and (\ref{neq2})
are immediate
implication of the bound (\ref{neq0})  and  the convergence of the $\beta'$s. 
Indeed,  
for any $\eps>0$ we have that 
$\sup_{x\in D}|\beta_{n(x,X(m))}(x)-\beta(x)|<\eps$
as soon as 
$\inf_{x\in D} n(x,X(m)) > n(\eps),$
for some $n(\eps)$.
Lemma \ref{L1} is proved.

{\bf Proof of Corollary \ref{c1}.} 
By the equation (\ref{f1}) the unconditional  
density of the random variable $X_{k+1}$ at point $x$ is 
$$\E \psi(x|X(k))=\E\frac{\beta_{n(x, X(k))}(x)}
{\alpha(X(k))}.$$
The integrand in this mean is 
bounded and converges in probability to $\beta(x)/\alpha$ as $k\ri \infty$
by Lemma \ref{L2}.
Therefore, $\E \psi(x|X(k))\ri \beta(x)/\alpha$ for any $x\in D$
 as $k\ri \infty$.
It is well known that the point-wise convergence of 
densities implies the convergence in total variation. 
Corollary \ref{c1} is proved.

{\bf Proof of Lemma  \ref{L2}.}
To simplify the  notation we assume that the
Lebesgue measure of the set $D$ is $1$.
We start with  part $1)$.
Let $\delta_{0}$ be a constant defined
in Lemma \ref{L1}. 
Note that 
$$J_k(f)=
\frac{1}{\alpha(X(k-1))} \int\limits_{D}f(x)\beta_{n(x,X(k-1))}(x)dx,\,
\,k\geq 1.$$
Fix an arbitrary $\eps>0$ 
and define
\begin{equation}
\label{event}
B_{k,\eps}=\left\{\sup\limits_{x\in D}|\beta(x)-\beta_{n(x,X(k-1))}(x)|\geq 
\eps \right\},\,\,\, k\geq 1.
\end{equation}
One can write
\begin{align*}
\E \left|J_k(f)-J(f)\right|^p&=\E \left|J_k(f)-J(f)\right|^p 
I_{\{B_{k,\eps}\}}+ \E \left|J_k(f) 
-J(f)\right|^p I_{\{\overline{B}_{k,\eps}\}}\\
&=S_1+S_2,
\end{align*}
where by $I_{\{B\}}$ we denoted an indicator of an event $B$. 
It is easy to see that
\begin{align}
\label{eks0}
J_k(f)-J(f) & =
\frac{\int_{D}f(x)(\beta_{n(x,X(k-1))}(x)-\beta(x))dx }{\alpha(X(k-1))}\\
&+J(f)\frac{\int_{D}(\beta(x)-\beta_{n(x,X(k-1))}(x))dx }{\alpha(X(k-1))}\nonumber,
\end{align}
hence
\begin{equation}
\label{eks}
|J_k(f)-J(f)|\leq  
\frac{2\|f\|_{\infty}}{\beta_{min}}
\sup\limits_{x\in D}|\beta(x)-\beta_{n(x,X(k-1))}(x)|.
\end{equation}
Let $k(\eps)$ be such that 
$\|\beta_k-\beta\|_{\infty}\leq \eps$ for any $k>k(\eps)$.
Then for any $k>k(\eps)$ we can bound
\begin{equation}
\label{eq3}
S_1\leq  C\eps^p.
\end{equation}
Using Lemma \ref{L1}  we have that for sufficiently large $k$
\begin{equation}
\label{eq4}
S_2\leq \left(\frac{4\|f\|_{\infty}\beta_{max}}{\beta_{min}}\right)^p
\P\{\overline{B}_{k,\eps}\}\leq Ce^{-\lambda k}.
\end{equation}
Combining bounds (\ref{eq3}) and (\ref{eq4}) we get that 
for all sufficiently large $k$ 
$$\|J_k(f) -J(f)\|_{L^p}^p\leq C(\eps^p+e^{-\lambda k}).$$
Therefore $L^p$-convergence of $J_k(f)$ to $J(f)$ is proved
for any $p>1$, since $\eps$ was taken arbitrary.
Part $1)$ of the lemma is proved.

Let now  the condition (\ref{var}) holds. 
Fix an arbitrary $\delta\in (0, \delta_0)$ and 
define  
\begin{equation}
\label{event1}
B_{k,\delta}=\left\{\inf\limits_{x\in D} n(x, X(k-1))\geq k\delta\right\},\,\,\,
k\geq 1.
\end{equation}
One can repeat the reasonings above using 
this sequence of events instead of the events (\ref{event}) and 
get the bound
$S_1\leq C\varphi^{p}(k\delta)$, therefore 
part $2)$ of Lemma \ref{L2} is also proved.

{\bf Remark.} Note  that in the equation (\ref{eks}) it is  
also possible  to bound
\begin{equation}
\label{eks00}
|J_k(f)-J(f)|\leq  
\frac{2\|f\|_{1}}{\beta_{min}}
\sup\limits_{x\in D}|\beta(x)-\beta_{n(x,X(k-1))}(x)|,
\end{equation}
where $\|f\|_{1}=\int_D|f(x)|dx$.

{\bf Proof of Corollary \ref{c2}.}
By the binomial formula we have that  
$$\left|\E\left(\xi_{k}^n(f)|{\cal F}_{k-1}\right)-{\cal U}_{n}(f)\right|
\leq \sum\limits_{i=0}^{n}{n\choose i} 
\left|J_{k}(f^i)(\E f(X_{k}))^{n-i}-J(f^i)J^{n-i}(f)\right|.$$
Noting that
\begin{align*}
\left|J_{k}(f^i)(\E f(X_{k}))^{n-i}-J(f^i)J^{n-i}(f)\right|
 & \leq 
C\left(\left|J_{k}(f^i)-J(f^i)\right|\right.\\
&+\left.\left|\E f(X_{k})-J(f)\right|\right)
\end{align*}
and applying part $1)$ of Lemma \ref{L2} we prove part $1)$ of the corollary.
If the condition (\ref{expon}) holds, then by part $2)$
of Lemma \ref{L2} we can bound for any $\delta\in (0,\delta_0)$
\begin{equation}
\E \left|J_{k}(f^i)-J(f^i)\right|
+\left|\E f(X_{k})-J(f)\right|\leq 
C\left(\varphi(k\delta)+e^{-\lambda k}\right)
\end{equation}
and part 2) of the corollary is also proved.

{\bf Proof of Lemma \ref{L3}.}
We can write 
\begin{align*}
\label{Cauc}
\E \left(\eta\prod\limits_{v=1}^{k}\xi^{r_v}_{i_v}(g_{v})\right) & =
{\cal U}_{r_k}(g_{k})
\E \left(\eta\prod\limits_{v=1}^{k-1}\xi^{r_v}_{i_v}(g_{v})\right)
 \\
& +\E \left(\eta \prod\limits_{v=1}^{k-1}\xi^{r_v}_{i_v}(g_{v})
\left(\E\left(\xi_{i_k}^{r_k}(g_k)|\F_{i_{k}-1}\right)- 
{\cal U}_{r_k}(g_{k})\right)\right).
\end{align*}
The functions $g'$s are bounded, so 
\begin{multline*}
\left| 
\E \eta \prod\limits_{v=1}^{k-1}\xi^{r_v}_{i_v}(g_{v})
\left(\E\left(\xi_{i_k}^{r_k}(g_k)|\F_{i_{k}-1}\right)- 
{\cal U}_{r_k}(g_{k})\right)
\right|\\
\leq C_{1}^{n-r_k} 
\E \left|\E\left(\xi_{i_k}^{r_k}(g_k)|\F_{i_{k}-1}\right)- 
{\cal U}_{r_k}(g_{k})\right|,
\end{multline*}
and the right hand side above goes to $0$ 
as $i_k\ri \infty$ by part $1)$ of Corollary \ref{c2}.
If the condition (\ref{var}) holds, then
by part $2)$ of Corollary \ref{c2} we can bound  
$$
\E \left|\E\left(\xi_{i_k}^{r_k}(g_k)|\F_{i_{k}-1}\right)- 
{\cal U}_{r_k}(g_{k})\right|
\leq C_2\left(\varphi(i_k\delta)+e^{-\lambda i_k}\right) 
$$ 
for any  $\delta\in (0,\delta_0)$ with some $\lambda=\lambda(\delta)$.
Repeating the same arguments for the indices $i_{k-1},\ldots,
i_1$ in 
$\E \left(\eta \prod_{v=1}^{k-1}\xi^{r_v}_{i_v}(g_{v})\right)$ 
we finish the proof.

{\bf Proof of Lemma \ref{L4}.}
Let us prove that 
\begin{equation}
\label{l41}
\frac{1}{\sqrt{m}}\sum\limits_{k=1}^{m}
\left(J_{k}(f)-J(f)\right) \ri 0,
\end{equation}
in probability as $m\ri \infty$.
Using the bound (\ref{eks})
we get that
$$
|J_k(f)-J(f)| \leq 
C_1\int\limits_{D}|\beta_{n(x,X(k-1))}(x)-\beta(x)|dx I_{\{B_{k,\delta}\}}
+C_2 I_{\{\overline{B}_{k,\delta}\}}
$$
where $B_{k,\delta}$ is the event defined by the equation (\ref{event1}). 
Therefore
\begin{align}
\frac{1}{\sqrt{m}}\left|\sum\limits_{k=1}^{m}
(J_{k}(f)-J(f))\right|& \leq 
\frac{C_1}{\sqrt{m}}\sum\limits_{k=1}^{m}\int\limits_{D}
|\beta_{n(x,X(k-1))}(x)-\beta(x)|dx I_{\{B_{k,\delta}\}} \label{l42}\\
& +
\frac{C_2}{\sqrt{m}}\sum\limits_{k=1}^{m}I_{\{\overline{B}_{k,\delta}\}}
\nonumber.
\end{align}
By Lemma \ref{L1}
$$\sum\limits_{k=1}^{\infty}\P\{\overline{B}_{k,\delta}\}<\infty,$$
hence by Borel-Cantelli lemma only a finite number
of events $\overline{B}_{k,\delta}$ occurs
with probability $1$, so 
$$\frac{C_2}{\sqrt{m}}\sum\limits_{k=1}^{m}I_{\{\overline{B}_{k,\delta}\}}
\ri 0$$
almost surely as $m\ri \infty$.
The first sum in the right hand side of the equation (\ref{l42})
is bounded by
\begin{displaymath}
\frac{C_1}{\sqrt{m}}\sum\limits_{k=1}^{m}
\sup\limits_{x\in D}|\beta_{n(x,X(k-1))}(x)-\beta(x)|I_{\{B_{k,\delta}\}}
 \leq
\frac{C_3}{\sqrt{m}}\sum\limits_{k=0}^{m}\varphi (k\delta),
\end{displaymath}
and it goes to $0$ as $m\ri \infty$ because of the 
equation (\ref{var}).
Repeating the same arguments we  can also prove that
$$\frac{1}{\sqrt{m}}\sum\limits_{k=1}^{m}
\left(\E f(X_k)-J(f)\right) \ri 0,$$
as $m\ri \infty$, therefore Lemma \ref{L4} is proved.

\section{Exponential rate of convergence}
If the rate of convergence in (\ref{expon}) is exponential, 
namely
if $\varphi(k)=\exp(-\gamma k)$ for some $\gamma>0$, then 
stronger statement of asymptotic independence
of random variables $X_k,\, k\geq 1,$ can be made.
Fix some $0<\eps<1/2$ and 
denote
$$\tilde{S}_m(f)=
\frac{1}{\sqrt{m-m^{\eps}}}
\sum\limits_{k=m^{\eps}}^{m}(f(X_k)-\E f(X_k)).$$
Let $Y_i,i\geq 1,$ be a collection of independent
random variables  with  the common probability density 
$\beta(x)/\alpha$.
Denote 
$$S_{0,m}(f)=\frac{1}{\sqrt{m_{\eps}}}
\sum\limits_{k=1}^{m_{\eps}}(f(Y_k)-\E f(Y_k)),$$
where we denoted $m_{\eps}=m-m^{\eps}$.
We are going to show that for a fixed set of positive indexes
$r_1,\ldots,r_k$, such that $r_1+\cdots+r_k=n$
the following expansion holds
\begin{equation}
\label{exp}
\prod\limits_{j=1}^{k}\E \tilde{S}_m^{r_j}(f)=
\prod\limits_{j=1}^{k}\E S_{0,m}^{r_j}(f)
+ \zeta_{m}(r_1,\ldots,r_k,f),
\end{equation}
where 
$$|\zeta_{m}(r_1,\ldots,r_k,f)|\leq C(n) m^{\eps+n/2} e^{-\rho m^{\eps}}.$$
For the simplicity of notation 
we prove the expansion (\ref{exp}) for 
the particular case $k=1, r_1=n$.
It is easy to see that 
$$\E \tilde{S}_m^{n}(f)=m^{-n/2}_{\eps}\sum\limits_{t_1,\ldots,t_p}
\,\,\, \sum\limits_{m_{\eps}\leq i_1<\ldots < i_p \leq m}
\E \prod\limits_{v=1}^{p}\xi^{t_v}_{i_v}(f),$$
where the first sum is over all sets of positive integers
$t_i,\, i=1,\ldots,p,$ such that $t_1+\cdots +t_p=n$. 
We get the expansion (\ref{exp}) if we put
 $$\zeta_{m}(n,f)=m^{-n/2}_{\eps}\sum\limits_{t_1,\ldots,t_p}
\,\,\, \sum\limits_{m_{\eps}\leq i_1<\ldots < i_p \leq m}
\left(\E \prod\limits_{v=1}^{p}\xi^{t_v}_{i_v}(f)-
 \prod\limits_{v=1}^{p}{\cal U}_{t_v}(f)\right).$$
Applying the bound (\ref{bb}) with $\varphi(k)=\exp(-\gamma k)$ yields 
that 
%\begin{equation}
%\label{bb1}
$$
\left|\E \prod\limits_{v=1}^{p}\xi^{t_v}_{i_v}(f)-
 \E \prod\limits_{v=1}^{p}{\cal U}_{t_v}(f)\right| 
\leq C\sum\limits_{v=1}^{p} e^{-\rho i_v},
$$
%\end{equation}
where $\rho=\min(\gamma,\lambda)$.
Therefore  we get that 
\begin{equation}
\label{zet} 
|\zeta_{m}(n,f)|\leq m^{-n/2}_{\eps}\sum\limits_{t_1,\ldots,t_p}
\,\,\, \sum\limits_{m_{\eps}\leq i_1<\ldots < i_p \leq m}
C\sum\limits_{v=1}^{p} e^{-\rho i_v}.
\end{equation}
It is easy to see that for any fixed set of positive integers 
$t_1,\ldots,t_p$ in the first sum we can bound
\begin{align*}
m^{-n/2}_{\eps}\sum\limits_{m_{\eps}\leq i_1<\ldots < i_p \leq m}
C_1^n\sum\limits_{v=1}^{p} e^{-\rho i_v} & \leq 
C_2 m^{(\eps-1/2)n}(m-m^{\eps})^{p-1}e^{-\rho m^{\eps}}\\
& \leq C_3 m^{\eps+n/2}e^{-\rho m^{\eps}}.
\end{align*}
The first sum
in (\ref{zet}) contains the number of terms depending only on $n$,
therefore
$$|\zeta_{m}(n,f)|\leq C_4(n) m^{\eps+n/2}e^{-\rho m^{\eps}}.$$
Using the representation (\ref{exp}) 
we can prove that ${\cal K}_{mn}(f)$  the $n$th cumulant 
of  $\tilde{S}_m(f)$ converges as $m\ri \infty$ 
to the cumulant of a Gaussian random variable with zero mean 
and the variance $G(f,f)$. 
Using Lemma \ref{L3} it is easy to prove  that 
${\cal K}_{m2}(f)\ri G(f,f)$ as $m\ri \infty$.
Let us to prove that  ${\cal K}_{mn}(f)\ri 0$ as
$m\ri \infty$ for $n\geq 3$. 
Recall that the cumulants ${\cal K}_{mn}(f),\, n\geq 1,$ 
are defined as the Taylor coefficients
of the logarithm of the characteristic function
\begin{equation}
\label{semi}
\log \E e^{it\tilde{S}_m(f)}=\sum\limits_{n=1}^{\infty}
{\cal K}_{mn}(f)\frac{(it)^n}{n!},\,\,\,t\in \R.
\end{equation}
%The series (\ref{semi}) absolutely converges,
%since the random variable $\tilde{S}_m(f)$ has all finite moments.
Each cumulant can be presented as a finite linear combination
of the products of moments (see, for instance, \cite{Leon})
 \begin{equation}
\label{semi1}
{\cal K}_{mn}(f)=\sum\limits_{k=1}^{n}(-1)^{k-1}(k-1)!
\sum\limits_{r_1,\ldots,r_k}\prod\limits_{j=1}^{k}\E \tilde{S}_m^{r_j}(f),
\end{equation}
where the second sum is over all sets of positive  integers 
$\{r_1,\ldots,r_k\}$ such that  
$r_1+\cdots+r_k=n$.
The equation (\ref{exp}) yields that
\begin{equation}
\label{exp2}
{\cal K}_{mn}(f)={\cal K}^{(0)}_{mn}(f)+
\sum\limits_{k=1}^{n}(-1)^{k-1}(k-1)!
\sum\limits_{r_1,\ldots,r_k} \zeta_{m}(r_1,\ldots,r_k,f),
\end{equation}
where 
${\cal K}^{(0)}_{mn}(f)$ is $n$th cumulant of the random variable
$S_{0,m}(f)$.
Because of the independence we have that
${\cal K}^{(0)}_{mn}(f)\sim m_{\eps}^{-n/2+1}\ri 0$
for any $n>2$ as  $m \ri \infty$.
It remains to note that
$$\left|\sum\limits_{k=1}^{n}(-1)^{k-1}(k-1)!
\sum\limits_{r_1,\ldots,r_k} \zeta_{m}(r_1,\ldots,r_k,f)\right|\leq
C(n)n!m^{\eps+n/2}e^{-\rho m^{\eps}} \ri 0,$$
as $m\ri \infty$. Thus the convergence of cumulants is proved.
It is well known that this implies the weak convergence. 
%(see, for instance,
%Lemma 3 in \cite{Sosh}).
\vspace{1mm}

\paragraph{{\large Acknowledgments.}}
This research was supported by the Technology Foundation 
STW, applied
science division of NWO, and the technology program of the Ministry of
Economic Affairs, The Netherlands
 (project CWI.6155 'Markov sequential point processes
for image analysis and statistical physics').

The author would like to thank Vadim Malyshev for motivating question
and Richard Gill for valuable comments and helpful  suggestions.

\end{document}